\documentclass[12pt]{article}
\usepackage[all]{xy}

\usepackage{epic}
\usepackage{amsmath}[1996/11/01]
\usepackage{amssymb,amsthm,amsfonts,latexsym,epsfig,graphics}
 \def\beql#1#2\eeql{\begin{equation}\label{#1}#2\end{equation}}

\advance \textheight by -1 \topmargin
\advance \textheight by -1 \headheight
\advance \textheight by -1 \headsep
\advance \textheight by -2 \footskip
\vsize = \textheight
\hsize = \textwidth
\DeclareMathOperator{\diag}{diag}

\DeclareMathOperator{\GL}{GL}
\DeclareMathOperator{\Aut}{Aut}
\DeclareMathOperator{\im}{im}

\DeclareMathOperator{\SL}{SL}
\DeclareMathOperator{\PSL}{PSL}

\newtheorem{theorem}{Theorem}[section]

\newtheorem{prop}[theorem]{Proposition}

\newcommand{\bew}{\noindent\underline{Proof.}\ }
\newtheorem{rem}[theorem]{Remark}
\newtheorem{lemma}[theorem]{Lemma}
\newtheorem{proposition}[theorem]{Proposition}
\newtheorem{definition}[theorem]{Definition}

\newcommand{\disj}{\stackrel{.}{\cup}}

\newcommand{\Z}{{\mathbb{Z}}}
\newcommand{\Q}{{\mathbb{Q}}}
\newcommand{\F}{{\mathbb{F}}}
\newcommand{\N}{{\mathbb{N}}}
\newcommand{\R}{{\mathbb{R}}}
\newcommand{\A}{{\mathbb{A}}}

\newcommand{\C}{{\mathbb{C}}}
\newcommand{\trace}{\mbox{trace}}
\newcommand{\tr}{\mbox{tr}}
\newcommand{\eb}{\phantom{zzz}\hfill{$\square $}\smallskip}

\renewcommand{\em}{\sf}

\begin{document}
\begin{center}
{\Large {\bf An even unimodular 72-dimensional lattice of minimum 8.}} \\
\vspace{1.5\baselineskip}
{\em Gabriele Nebe} \\
\vspace*{1\baselineskip}
Lehrstuhl D f\"ur Mathematik, RWTH Aachen University\\
52056 Aachen, Germany \\
 nebe@math.rwth-aachen.de \\
\vspace{1.5\baselineskip}
\end{center}

{\sc Abstract.}
{\small 
An even unimodular 72-dimensional lattice $\Gamma $ having minimum 8 is 
constructed as a tensor product of the Barnes lattice and 
 the Leech lattice over the ring of integers in the 
imaginary quadratic number field with discriminant $-7$.
The automorphism group of $\Gamma $ contains the 
absolutely irreducible rational matrix group 
$(\PSL_2(7) \times \SL _2(25)) : 2$.
}

\section{Introduction.}

In this paper a lattice $(L,Q)$ is always an even positive definite lattice,
i.e. a free $\Z $-module $L$ equipped with a quadratic form $Q:L\to \Z $
such that the bilinear form 
$$(\cdot,\cdot ): L\times L \to \Z , (x,y):=Q(x+y)-Q(x)-Q(y) $$
is positive definite on the real space $\R\otimes L$.
The dual lattice is
$$L^{\#} := \{ x\in \R \otimes L \mid (x,\ell) \in \Z \mbox{ for all } \ell \in L \} $$ 
and $L$ is called {\em unimodular}, if $L=L^{\#}$. 
The {\em minimum} of $L$ is  twice the minimum of the quadratic 
form on the non-zero vectors of $L$ 
$$\min (L) =  \min \{ (\ell ,\ell ) \mid 
0\neq \ell \in L \} .$$
From the theory of modular forms it is known (\cite{Siegel},   \cite{MallowsSloane})
that the minimum of an even unimodular lattice of dimension $n$ 
is always $\leq 2 \lfloor \frac{n}{24} \rfloor + 2  $. 
Lattices achieving this bound are called {\em extremal}. 
Of particular interest are extremal unimodular lattices in the 
so called ``jump dimensions'', these are the multiples of 24.
There are four even unimodular lattices known in the jump dimensions,
the {\em Leech lattice} $\Lambda $,
the unique even unimodular 
lattice in dimension 24 without roots, and three lattices called 
$P_{48p}$, 
$P_{48q}$, 
$P_{48n}$,  of dimension 48 which have minimum 6 \cite{ConwaySloane}, \cite{cycloquat}. 

It was a long standing open problem whether there exists 
an extremal 72-dimensional unimodular lattice 
(\cite[p. 151]{Sloane}, \cite[Section 3.4]{SchaSchuPi}).
Many people tried to construct such a lattice, 
 or to prove its non-existence.
Most of these attempts are not documented, all constructed lattices
contained vectors of norm 6. 
In 
\cite{BachocNebe} 
Christine Bachoc and I 
 discovered two extremal lattices in dimension 80 
of which we could prove extremality
using a classical construction that we learned from \cite{Queb}.
Given a binary code $C\leq \F_2^d$ and an 
even lattice $L\leq \R^n$ of odd determinant together with a polarisation 
$L/2L = T_1\oplus T_2$ by isotropic subspaces 
the new lattice ${\cal L}(T_1,T_2,C)$ of dimension $dn$ is constructed as 
the preimage in $L^d$ of 
$T_1\otimes C \oplus T_2 \otimes C^{\perp } \leq (L/2L)^d $. 
Inspired by \cite{GrossElkies} 
Christine and I used polarisations coming from  Hermitian 
$\Z [\alpha ]$-structures (where $\alpha ^2-\alpha +2 = 0$)
of $L$. 

Bob Griess' article \cite{Griess} 
%\footenote{I reviewed it for the Zentralblatt}
analyses this construction for certain polarisations 
of the Leech lattice $\Lambda $ and $C=\langle (1,1,1) \rangle \leq \F_2^3$ 
for which  he describes 
a strategy to prove extremality of the resulting lattice. 
This motivated me to try the
nine $\Z[\alpha ]$-structures of  $\Lambda $
calculated in \cite{Hentschel}.
I computed the number of vectors of norm 6 in all nine 72-dimensional lattices
using four different strategies: 
A combination of lattice reduction programs applied directly 
to the 72-dimensional lattice 
found vectors of norm 6 for all but one lattice.
I then went on to compute the super offenders as described in
\cite[Section 4]{Griess} and computed the number of norm 6 vectors
in the lattices as given in Table 1.
Using the $\Z[\alpha ]$ structure of $\Lambda $ this computation
may be reduced to a computation within the set of minimal
vectors of the Leech lattice.
 %as described in a first version of this paper.
The result of these computations agreed with the ones applying the methods
given in Section \ref{Minimum}.

Using the explicit matrices for this extremal lattice $\Gamma $  and the 
action of the subgroup $G$ of $\Aut (\Gamma )$ 
as constructed in Section \ref{Herm} 
Mark Watkins (personal communication) succeeded in listing
representatives of all $G$-orbits of the vectors of norm 8 in $\Gamma $
using the method described in \cite{StehleWatkins}.
From the stored information one verifies that 
$\Gamma $ has $6,218,175,600$  minimal vectors 
 which gives an
independent proof of the extremality of $\Gamma $ (see also
Theorem \ref{norm68} for an explicit description of the kissing
configuration of $\Gamma $).

\section{An Hermitian tensor product construction of $\Gamma $.}\label{Herm}

Throughout the paper let $\alpha $ be a generator of 
the ring of integers $\Z [\alpha ]$
in the imaginary quadratic number field  of discriminant $-7$,
with $\alpha ^2-\alpha +2 = 0$ and $\beta := \overline{\alpha } = 1-\alpha $ 
its 
complex conjugate.
Then $\Z [\alpha ] $ is a principal ideal domain and 
$(\alpha )$ and $(\beta )$ are the two maximal ideals 
of $\Z[\alpha ]$ that contain 2.

Let $(P,h)$ be an Hermitian $\Z [\alpha ]$-lattice,
so $P$ is a free $\Z[\alpha ]$-module  and 
$h:P\times P \to \Z [\alpha ] $ a positive definite Hermitian form. 
One example of such a lattice is the
 {\em Barnes lattice} 
$$P_b = \langle (\beta ,\beta ,0), (0,\beta ,\beta ), (\alpha,\alpha,\alpha )
\rangle = \langle (1,1,\alpha ),(0,\beta,\beta), (0,0,2 ) \rangle 
\leq \Z[\alpha ]^3 $$ 
with the half the standard Hermitian form 
$$h:P_b\times P_b \to \Z[\alpha ], h((a_1,a_2,a_3),(b_1,b_2,b_3)) = 
\frac{1}{2} \sum _{i=1}^3 a_i \overline{b_i} .$$
Then $P_b$ is Hermitian unimodular, 
$P_b = P_b^* := \{ v\in \Q P_b \mid h(v,\ell ) \in \Z[\alpha ] \mbox{ for all }
\ell \in P_b \} $.
The automorphism group of the $\Z[\alpha ]$-lattice $P_b$ is 
isomorphic to $\pm \PSL_2(7)$ 
$$\Aut _{\Z[\alpha ]} (P_b) := \{ g\in \GL (P_b) \mid 
h(gv,gw) = h(v,w) \mbox{ for all } v,w \in P_b \}  \cong \pm \PSL _2(7) .$$

From any such Hermitian $\Z[\alpha ]$-lattice $(P,h)$
one obtains an even $\Z $-lattice 
$$L(P,h) := (L,(,)) := (P,\trace _{\Z[\alpha ] / \Z } \circ h)$$
by restricting scalars. 
Since $h$ is Hermitian $h(\ell ,\ell) \in \Z = \Z[\alpha ] \cap \R $
 for all $\ell \in P$ and hence 
$$Q(\ell) := \frac{1}{2} \trace _{\Z[\alpha ] / \Z } (h(\ell,\ell)) 
= h(\ell ,\ell) \in \Z .$$
The dual lattice of $L(P,h)$ is the product of $P^*$ with the 
different of $\Z[\alpha ]$:
$$L(P,h)^{\#} := \{ v \in \Q P \mid \trace _{\Z[\alpha ] / \Z } (h(\ell,v)) 
\in \Z \mbox{ for all } \ell \in P \} = \frac{1}{\sqrt{-7}} P^* .$$

Michael Hentschel
\cite{Hentschel} classified all Hermitian $\Z[\alpha ]$-structures 
on the 
even unimodular $\Z $-lattices of dimension 24 using the 
Kneser neighbouring method (\cite{Kneser}) to generate the lattices 
and checking completeness with the mass formula.
In particular there are exactly nine such $\Z[\alpha ]$ structures 
$(P_i,h)$  ($1\leq i \leq 9$)
such that 
$(P_i, \frac{1}{7} \trace _{\Z[\alpha ]/\Z } \circ h) \cong \Lambda $ is
the Leech lattice.
The 36-dimensional 
Hermitian $\Z[\alpha ]$-lattice $R_i $ is defined as
$$(R_i,h) := P_b \otimes _{\Z[\alpha ]} P_i , \mbox{ so  } 
\Aut _{\Z[\alpha ]}  (R_i) \supseteq \PSL_2(7) \times \Aut _{\Z[\alpha ]} (P_i)  .$$

\begin{definition}
For $1\leq i\leq 9$ let 
$(\Gamma _i,(,)) := L(R_i,\frac{1}{7} h) := (R_i, \frac{1}{7} \trace _{\Z[\alpha ]/\Z } \circ h) $
where the quadratic form is $Q(\ell ) =\frac{1}{14}\trace_{\Z[\alpha ]/\Z }
 (h(\ell,\ell)) $ for all $\ell \in \Gamma _i = R_i$.
\end{definition}

All $\Gamma _i$ are even unimodular lattices of dimension 72. 
%In fact, the lattices $\Gamma _i$ may be constructed from 
%a polarisation of the Leech lattice as defined by Bob Griess
%\cite{Griess}. This construction is repeated in Section \ref{SecGriess}.
%It shows that $\min (\Gamma _i) \geq 6$. 
%We use Griess' construction to compute the number of norm 6 vectors in $\Gamma _i$. 

The table below lists these nine Hermitian structures of the
Leech lattice. The first column gives the structure of 
the automorphism group $\Aut _{\Z[\alpha ] }(P_i )$ 
followed by its order and then the number of vectors of norm 6
in the lattice $\Gamma _i$ (computed in Section \ref{Minimum} below).
\begin{center} 
Table 1 \\
$
\begin{array}{|c|c|c|c|} 
\hline 
 & \mbox{ group } & \mbox{ order } & \mbox{ norm 6 vectors } \\
\hline 
1 & \SL _2(25) &  2^4  3 \cdot 5^2  13 & 0 \\
\hline
2 & 2.A_6\times D_8 & 2^7  3^2  5 & 2\cdot 20,160 \\ 
\hline
3 & \SL_2(13).2 &  2^4  3 \cdot 7 \cdot 13 & 2\cdot 52,416 \\ 
\hline
4 & (\SL_2(5)\times A_5).2 &  2^6  3^2  5^2 & 2\cdot 100,800 \\
\hline
5 & (\SL_2(5)\times A_5).2 &  2^6  3^2  5^2 & 2\cdot 100,800 \\
\hline
6 & \mbox{ soluble } & 2^93^3 & 2\cdot 177,408 \\
\hline
7 & \pm \PSL_2(7) \times  (C_7:C_3) &  2^4 3^2  7^2 & 2 \cdot 306,432 \\ 
\hline
8 & \PSL _2(7) \times  2.A_7 &  2^7  3^3  5 \cdot 7^2  & 2\cdot 504,000 \\
\hline
9 & 2.J_2.2 & 2^9  3^3  5^2 7  & 2 \cdot 1,209,600 \\
\hline
\end{array}
$
\end{center}

\begin{rem} 
(a) 
The groups number $1$, $3$, $4,5$, and $9$ are maximal finite quaternionic
matrix groups with endomorphism algebra the definite quaternion algebra
with center $\Q $ and discriminant $5^2$ ($1,4,5,9$) resp. 
$13^2$ (group number $3$) (see \cite{quat}).
For the group number $4$ resp. $5$, the endomorphism ring of the
lattice is not the maximal order. \\
(b) The group number $8$  is a maximal finite 
symplectic matrix group over $\Q[\alpha ]$ as defined in \cite{Kirschmer},
it is globally irreducible in the sense of \cite{Gross}.  \\
(c) The groups number $2$ and $7$ are reducible. 
\end{rem} 

The Hermitian structures number 4 and 5 are 
just Galois conjugate to each other, whereas all the others 
are Galois invariant. 
For these seven lattices 
the automorphism group of the $\Z $-lattice $\Gamma _i$ hence contains an
extension of $\Aut_{\Z[\alpha ]} (R_i)$ by the Galois automorphism. 
For the extremal lattice $\Gamma :=\Gamma _1$ this is a split extension.

\begin{theorem}
The lattice $\Gamma $ is an extremal even unimodular 
lattice of dimension $72$. 
Its automorphism group 
$\Aut (\Gamma )$ contains the subgroup ${\cal U}:= (\PSL_2(7) \times \SL_2(25) ) :2  $.
\end{theorem}

Two  proofs that the minimum of $\Gamma $ is 8
are given below.

\begin{rem}
The natural $\C {\cal U} $-module $\C \otimes \Gamma $ contains 
no ${\cal U}$-invariant submodules, so 
$\Aut (\Gamma )$ is an absolutely irreducible 
subgroup of $\GL_{72} (\Q )$. 
In fact ${\cal U}$ is almost a globally irreducible representation
in the sense of \cite{Gross}. 
More precisely $\F_p \otimes \Gamma $ is also absolutely irreducible 
except for $p=5$ and $p=7$, where the module has a unique non-trivial
submodule, which is of dimension $36$. 
For both primes $p=5$ and $p=7$ there is an element $x_p \in 
N_{\GL_{72}(\Q )} ({\cal U})$, the rational normalizer of ${\cal U}$
mapping $\Gamma $ to the unique sublattice of index $p^{36}$,
which is therefore isometric to $(\Gamma , pQ)$ (see
\cite{normaliser}).
Therefore $\Aut (\Gamma )$ %(which very likely equals ${\cal U}$)
is a maximal finite subgroup of $\GL_{72}(\Q )$.
\end{rem} 

\begin{rem}
Since $\Aut (\Gamma )$ contains an element of order $91$ the 
lattice $\Gamma $ is an ideal lattice in the cyclotomic
field $\Q [\exp (2\pi i/91 ) ] $ in the sense of \cite{Bayer}.
It would be interesting  to determine the ideal class 
of this lattice.
\end{rem}

\subsection{An elementary linear algebra construction.}

This section just repeats the construction above in 
elementary linear algebra (understood by computer algebra
systems). 

Let $(b_1,\ldots , b_{24})$ be a $\Z $-basis of the Leech lattice
$\Lambda $
and $F:=((b_i,b_j)) \in \Z ^{24 \times 24} _{sym}$ denote its
Gram matrix. 
Then an Hermitian structure over $\Z [\alpha ]$ is given by a 
matrix $A\in \Z^{24 \times 24} $ such that 
$AFA^{tr} = 2 F$ and the $F$-adjoint $FA^{tr} F^{-1} = 1-A =: B$.
Mapping $\alpha $ to the right multiplication by $A$ then 
defines the Hermitian $\Z[\alpha ]$-structure on the $\Z $-lattice
$\Lambda $.

That there are exactly nine such $\Z[\alpha ]$ structures 
of the Leech lattice  means that there are nine such matrices
$A_1,\ldots , A_9$ up to conjugation under the automorphism group 
of $\Lambda $.

For any of these nine structures 
the even unimodular lattice $\Gamma _i $ of 
dimension 72 is constructed 
as a sublattice of $\Lambda \perp \Lambda \perp \Lambda $
with Gram matrix $\frac{1}{2} \diag (F,F,F)$ generated by the 
rows of the block matrix 
$$ \left( \begin{array}{ccc} A_i & A_i & A_i \\
B_i & B_i & 0 \\ 0 & B_i & B_i \end{array} \right) 
\mbox{ or equivalently } 
 \left( \begin{array}{ccc} 1 & 1 & A_i \\
0 & B_i & B_i \\ 0 & 0 &2 \end{array} \right)  =: T_i
$$

If $U$ denotes the subgroup of $\GL_{72} (\Z )$ 
obtained by replacing $\alpha $ by $A_i$ in the 
group 
$\Aut _{\Z[\alpha ]} (P_b) \leq \GL_3(\Z[\alpha ]) $ 
isomorphic to $\PSL_2(7) $
then $\Aut (\Gamma _i) $ contains the matrix group 
$$ \langle \{ \diag (g,g,g) \mid g\in \Aut(\Lambda ), g A_i = A_i g \} 
\cup U \rangle \cong \Aut _{\Z[\alpha ]}(P_i ) \times \PSL_2(7) .$$

A matrix for the additional Galois automorphism with respect to 
the basis given by $T_i$ above can be constructed from 
an isometry 
$$ Y_i \in \GL_{24}(\Z),\ \  Y_i F Y_i^{\tr} = F, \ \ Y_i A_i Y_i^{-1} = B_i $$ 
(this only exists for $i\neq 4,5$) 
as the block matrix 
$$ \left( \begin{array}{ccc} 
Y_i & -Y_i & A_iY_i  \\ 
0 & -B_i Y_i & Y_i \\
-AY_i & 0 & Y_i \end{array} \right) .
$$
The shape of the matrix was obtained from an 
isometry between $(P_b,h)$ and $(P_b,\overline{h}) $.

\subsection{A classical coding theory construction.} \label{codesconst} 

The lattices $\Gamma _i$ can be obtained using a special case
of a classical construction with codes: 
%In \cite{Griess} Bob Griess gives a construction of even unimodular
%lattices using  polarisations that also appeared in \cite{Queb}.
If $(L,Q)$ is an even unimodular lattice, then $L/2L$ becomes a
non-degenerate quadratic space over $\F_2$ with 
quadratic form $q(\ell +2L) := Q(\ell ) + 2\Z .$
This has Witt defect 0, so there are totally isotropic subspaces
$U,V \leq (L/2L,q) $ such that $L/2L = U\oplus V$. 
Let $2L\leq M,N \leq L$ denote the preimages of $U,V$, respectively.
Then $(M,\frac{1}{2} Q)$ and $(N,\frac{1}{2} Q)$ are again 
even unimodular lattices.

\begin{definition} (\cite{Griess}, \cite[Construction I]{Queb})
Given such a polarisation $(M,N)$ of the even unimodular 
lattice $(L,Q)$ and some $k\in \N$ let
$$L(M,N,k) := \{ (x_1+y,x_2+y,\ldots , x_k+y ) \mid 
y\in N, x_1,\ldots ,x_k\in M \mbox{ and } x_1+\ldots +x_k \in M\cap N \} .$$
Then the lattice $(L(M,N,k),\tilde{Q} )$ is an even unimodular lattice
(\cite[Proposition]{Queb})
where $$\tilde{Q} (x_1+y,x_2+y,\ldots , x_k+y ) := 
\frac{1}{2} \sum _{i=1}^k Q(x_i+y ) .$$
\end{definition}

From the explicit basis of the Barnes lattice given in Section \ref{Herm} 
one immediately sees the following.

\begin{rem}
Assume that the lattice $L$ has an Hermitian structure over 
$\Z[\alpha ]$ as defined in Section \ref{Herm}. 
Then $M:=\alpha L$ and $N:=\beta L$ defines a polarisation of $L$ such
that $L(M,N,3) \cong P_b \otimes _{\Z[\alpha ]} L $.
\end{rem}

\begin{rem} \label{min6} (\cite[Theorem 4.10]{Griess})
Assume that $L = \Lambda \cong (M,\frac{1}{2} Q) \cong (N,\frac{1}{2} Q)$
is the Leech lattice. 
Then $L(M,N,3) $ has minimum 6 or 8.
The vectors of norm 6 in $L(M,N,3) $ 
are of the form 
$(w+x,w+y,w+z)$ with $w\in N$, $x,y,z \in M$, $x+y+z \in 2\Lambda $ 
and $Q(w+x) = Q(w+y) = Q(w+z) = 2$.
\end{rem} 

The next section analyses this construction in particular for this
situation. Using properties of the Leech lattice a description of the
vectors of norm 6 and 8 in $L(M,N,3)$, 
in particular of the kissing configuration of the 
new extremal lattice, is obtained. 
I thank Noam Elkies for proposing this interesting 
question.

\section{The vectors of norm 6 and 8 in $L(M,N , 3)$.}\label{NORM8}

Let $(\Lambda ,Q) $ be the Leech lattice, 
and fix a polarisation $(M,N) $ of $\Lambda $ 
such that $(M,\frac{1}{2} Q) \cong (N,\frac{1}{2} Q) \cong \Lambda $. 

The following property of the Leech lattice is well known.

\begin{lemma}  \label{classes}
The nontrivial classes of $\Lambda / 2\Lambda $ are 
represented by vectors $v\in \Lambda $ of norm $(v,v) = 4,6 $ and $8$.
In particular all classes of $M/2\Lambda $ and $N/2\Lambda $ 
are represented by vectors of norm $8$.
If $K = v + 2\Lambda $ contains a vector of norm $8$, then 
$ \{ k\in K \mid (k,k) \leq 8 \} = \{ \pm k_1,\ldots , \pm k_{24} \}$ 
with $(k_i,k_j) = 8\delta _{ij}$. 
If $K=v+2\Lambda $ contains a vector $v$ of norm $4$ or $6$ then
$ \{ k\in K \mid (k,k) \leq 8 \} = \{ \pm v \}$.
\end{lemma}

\bew
Let $(v,v) = 4$ or $(v,v)=6$ and assume that there is some 
$\pm v \neq k\in v+2\Lambda $ such that $(k,k)\leq 8$. 
Then one of $v\pm k \in 2\Lambda $ has norm $\leq 6+8 < 16 = \min (2\Lambda )$ 
which is a contradiction.
Similarly one sees that for $(v,v)=8$ the vectors of norm $8$ in
$v+2\Lambda $ form a frame. 
Now 
$$\frac{|\Lambda _4|}{2} + 
\frac{|\Lambda _6|}{2} + 
\frac{|\Lambda _8|}{48}   = 2^{24}-1 $$ 
so all nonzero classes of $\Lambda /2\Lambda $ are represented by vectors of 
norm $\leq 8$.
\eb

\begin{proposition}\label{w3}
Fix some $w\not\in M$ with $(w,w)=8$. Then
$W_3(w) := \{ x\in M \mid Q(x+w) = 3 \} $
has cardinality
$2\cdot 2048 $.
\end{proposition}

\bew
If $x\in M$ such that $Q(w+x) = 3$ then the class $K_x:=x+2\Lambda $ 
is not perpendicular to $K_w$. Since $w\not\in M$ there are 
$2^{11} = 2048$ classes $x+2\Lambda \in M/2\Lambda $  that are not perpendicular
to $K_w \in N/2\Lambda $ by the polarisation property. 
So there are $2048$ possibilities for the class $K_{w+x} = K_x+K_w = 
(x+w) + 2\Lambda $. 
Such a class is necessarily anisotropic and therefore contains 
exactly 2 vectors of norm 6 by Lemma \ref{classes}.
\eb

\begin{rem}\label{countst} 
Let $\tilde{X} := \{ m\in M
\mid (m,m) = 8 \}$ denote the set of minimal vectors in $M$.
Then $\tilde{X} $ is a spherical 11-design. 
For any $w\in N$ with $(w,w) = 8$ and $i\in \Z _{\geq 0} $ let
$$n_i := |\{ x\in \tilde{X} \mid (x,w) = \pm i \} |$$ 
Then $n_i=0$ if $i\geq 7$,
$$\sum _{i \mbox{ odd }} n_i = 2048\cdot 48 \mbox{ and }
\sum _{i \mbox{ even }} n_i = 2047\cdot 48  .$$
This allows to compute all $n_i$ using the 11-design properties of $\tilde{X}$
(see for instance \cite{VenkovEnsMath} for a description of this method):
$$\begin{array}{|c|c|c|c|c|c|c|c|}
\hline
i & 0 & 1 & 2 & 3 & 4 & 5 & 6  \\
\hline
n_i & 
46488 & 78848 & 47216 & 18944 & 4536 & 512 & 16 \\
\hline
\end{array}
$$
In particular $S(w):= \{ x\in \tilde{X} \mid (x,w) = -6 \} $
contains exactly 8 elements. 
The only property we need is that this set is not empty.
\end{rem}

\begin{proposition}\label{w2} 
Fix some $w \in N $ with $(w,w)=8$. Then
$W_2(w) := \{ x\in M \mid Q(x+w) = 2 \} $
has cardinality
$2\cdot 24$.
The set $\{ x+w \mid x\in W_2(w) \} $ is the rescaled root system 
$24 \A _1$. 
\end{proposition}

\bew
Since $w\not\in M$ and $M/2\Lambda $ is maximal isotropic there
is some $m\in M$ such that $(m,w) $ is odd and in particular 
$M_w := \{ m\in M \mid (m,w) \mbox{ even } \} $ is a sublattice
of index 2 in $M$ and 
$M^w := \langle M_w, w \rangle = M_w \disj w + M_w $ is a neighbor 
of $M$  in the sense of \cite{Kneser} 
and again isometric to a rescaled even unimodular lattice.
Since  $M$ has no roots 
 the root system of $M^w$ is either empty or 
$24 \A _1 $ (rescaled to have minimum 4). 
Remark \ref{countst} shows that this lattice contains at least 
16 roots, so there are 
24 pairs of orthogonal vectors of norm 4 in $M^w$  giving rise to 
the $x+w$ with  $x\in W_2(w)$.
\eb

\begin{theorem}\label{norm68}
The vectors of norm 6 and 8 in $L(M,N,3)$ are of the form 
$(w+x,w+y,w+z)$ with $w\in M$, $x,y,z\in N$, $x+y+z \in 2\Lambda $ such that 
$Q(w+x) + Q(w+y) +Q(w+z) = 6$ or $8$. 
For norm 6 the only possibility is 
$Q(w+x) = Q(w+y) =Q(w+z) = 2$. Let $b_6$ denote the number of such vectors. 
For norm 8 one has the possible types 
\begin{itemize} 
\item[(a)] $(8,0,0)$ with $196560 \cdot 3$ vectors. 
\item[(b)] $(4,4,0)$ with $196560 \cdot 48 \cdot 3$ vectors. 
\item[(c)] $(3,3,2)$ with $4095 \cdot 48 \cdot 2048 \cdot 2\cdot 2 \cdot 3$ vectors. 
\item[(d)] $(4,2,2)$ with $4095 \cdot 48 \cdot 48 \cdot 48\cdot  3-72b_6$ vectors. 
\end{itemize}
\end{theorem}

\bew
Clearly these are the only possibilities for vectors of norm 6 or 8
in $L(M,N,3)$. 
The vectors of type $(8,0,0)$ correspond to the minimal vectors 
in the sublattice $2\Lambda \perp 2\Lambda \perp 2\Lambda $. 
The vectors of  type $(4,4,0)$ are of the form 
$(x,y,0)$ with $x,y\in M$, $(x,x)=(y,y) = 8$
such that $x+2\Lambda = y + 2\Lambda $, 
so one has $196560$ possibilities for $x$ and for each such $x$ one
may choose all $48$ minimal vectors $y \in x+2\Lambda $. 
The additional factor 3 counts the possible permutations
$(x,y,0), (x,0,y)$, and $(0,x,y)$. 
\\
To see (c) we first note that by Lemma \ref{classes} all 
anisotropic classes in $\Lambda /2\Lambda $ are represented 
by vectors of norm 6. 
For a fixed representative $w $ of one of the $4095$ classes of 
$N/2\Lambda $ of norm $(w,w) = 8$ we have to run through all 
$z\in W_2(w)$ and $y\in W_3(w)$. 
Then the condition that $x+y+z \in 2\Lambda $ means that 
$w+x+2\Lambda = y+z + 2\Lambda $, so $w+x$ is one of the 2 
vectors of norm 6 in this anisotropic class. 
\\
The last set is the set of vectors of Type $(4,2,2)$. 
Again, fixing some $w$ as above the elements 
$y$ and $z$ are in $W_2(w)$. So we have 48 possibilities of each of them.
Then $w+x \in w+y+z + 2\Lambda $ is in an isotropic class of
 $\Lambda/ 2\Lambda $.
This class is either of minimum 4 and then we have 2 vectors in this
class resulting in a vector of norm 6 in $L(M,N,3)$ or this
class is of minimum 8 and there are 48 possibilities for $w+x$ such that 
$Q(w+x) =4$. 
\eb

\section{Two proofs that the minimum of $\Gamma $ is 8} \label{Minimum} 

%To show that the minimum of $\Gamma $ is at least 8, 
%I used four different strategies. 
%First I applied a combination of lattice reduction programs to 
%all nine 72 dimensional lattices and found vectors of norm 6 
%for all but one lattice. 
%I then went on to compute the super offenders as described in 
%\cite[Section 4]{Griess} and computed the number of norm 6 vectors
%in the lattices as given in Table 1. 
%Using the $\Z[\alpha ]$ structure of $\Lambda $ this computation
%may be reduced to a computation within the set of minimal 
%vectors of the Leech lattice as described in a first version of 
%this paper. 
%I doublechecked all these computations by applying the method 
%given in Section \ref{Seoul} as well as the explicit description of
%vectors of norm 6 and 8 given in Section \ref{NORM8}. 
%

\subsection{Counting the norm 6 vectors in $L(M,N,3)$} \label{norm6} 

Let $W:= \{ w_1,\ldots ,w_{4095} \} $ denote a fixed set of 
representatives of the classes in $N/2\Lambda $ consisting of 
vectors of norm 8. 
The vectors of norm 6 in $\Gamma $ are of the 
form $(w+x,w+y,w+z)$ where $w \in W$, $x,y,z \in W_2(w) $ and
$x+y+z \in 2\Lambda $. 
The set $W_2(w)$ may be computed as  the set 
$W_2(w) = \{ s-w \mid s \in M^w , (s,s) = 4 \} $ 
as described in Proposition \ref{w2}. 

\begin{rem}
To count the vectors of norm 6 and 8 in $\Gamma _i$  let 
 $w\in W$ run through representatives of 
the $\Aut _{\Z[\alpha ]} (P_i)$-orbits on $N/2\Lambda $.
\\
For each such $w$ compute the set 
$W_2(w) = \{ s-w \mid s \in M^w , (s,s) = 4 \} $. 
Then run through the pairs 
$(x,y) \in W_2(w)^2$ and compute the vectors of norm 
$4$ in the lattice $\langle 2\Lambda , w+x+y \rangle $. 
This lattice either has two vectors of norm $4$ contributing to 
the vectors of norm $6$ in $\Gamma _i$ or 
it has minimum $8$ and then it contains $48$ vectors of norm $8$ 
contributing to the vectors of norm $8$ in $\Gamma _i$.
\\
Using this method I found the number of vectors of norm $6$ in
$\Gamma _i$ as given in Table 1.
%In particular all such lattices $\Gamma _i$ contain 
%at most  $4095\cdot 48\cdot 48\cdot 2 = 18,869,760$ vectors of
%norm 6.
%(this bound can be improved, if $x+y+z \in 2\Lambda $,
%then one of them is in $S(w) \cup -2w-S(w) $ and the
%other two come from $T(t)$, this yields 
%$4095\cdot 16\cdot 32\cdot 6 = 12,579,840$ for the maximal 
%cardinality).
\end{rem} 

\begin{rem}
As proposed  by Bob Griess in \cite[Lemma B.3]{Griess} 
 I checked all pairs of fourvolutions
 $f,g \in \Aut(\Lambda )$
such that $x:=fg$ is an element of odd prime order $p$ with irreducible
minimal polynomial. 
Then  $p=3,5,7,13 $ and the conjugacy class of $x$ in 
$2.Co_1$ is unique. 
To enumerate all such pairs, I computed the normaliser $N$ of $\langle x \rangle$
in $2.Co_1$ and went through all conjugacy classes of elements $f$ of $N$ 
such that $f^2=-1$. For each such $f$ that satisfies $x^f=x^{-1}$ 
I put $g:=f^{-1}x$. The centraliser in $2.Co_1$ acts on the situation.
All six lattices $L(\Lambda (f-1),\Lambda (g-1),3) $ contain
vectors of norm 6.
\end{rem}

%I also checked g=f^{-1}x^k for k \leq (p-1)/2 and interchanged f and g. 
%I always got the same number of norm 6 vectors, though in different C-orbits

\begin{center} 
Table 2 \\
$
\begin{array}{|c|c|c|} 
\hline 
  x & \mbox{ centraliser order } & \mbox{ norm 6 vectors } \\
\hline 
3 & 2^83^35^27 &   2\cdot 1,209,600 \\
5a & 2^83^35^27 &     2\cdot 1,209,600 \\
5b & 2^83 \cdot 5 &  2\cdot 103680 \\
7 & 2^33 \cdot 5 &  2\cdot 11520 \\
13 & 2^3 3 &   2\cdot 57600 \\
\hline
\end{array}
$
\end{center}

\subsection{Using orthogonal decomposition $(24,48)$.} \label{Seoul}

The idea (see \cite{seoul})
 is to embed the lattice $\Gamma _i$ into an
orthogonally decomposable lattice $I_1\perp I_2$ such that the
minimal vectors of $I_1$ and the minimum of sublattices of
$I_2$ can be computed.
The even unimodular
lattice $\Gamma _i$ has basis matrix $T_i$ given in Section \ref{Herm}
with respect to the Gram matrix $\frac{1}{2} \diag (F,F,F) $.
Let $\pi $ denote the orthogonal projection onto the first 24 components,
$K_1:=\ker(1-\pi) \subset I_1 := \im (\pi ) = K_1^{\#}  $, and 
$K_2:=\ker(\pi) \subset I_2 := \im (1-\pi ) = K_2^{\#} $.
Then $$ K_1 \perp K_2 \subset \Gamma _i  \subset  I_1 \perp I_2.$$
So the even unimodular lattice $\Gamma _i$ contains the sublattice
$K_1 \perp K_2$ of index $2^{24}$ and is contained in
$I_1 \perp I_2$ also of index $2^{24}$.
Moreover
$I_1 \cong \frac{1}{\sqrt{2}} \Lambda $ is isometric to
the Leech lattice scaled to have minimum 2  and
$I_2$ is a non integral lattice of dimension 48, with $(\ell ,\ell ) \in \Z $
for all $\ell \in I_2$.
A computer calculation shows that $\min (I_2) = 4$,
$K_1 \cong \sqrt{2} \Lambda $ has minimum 8 and
 also the 48-dimensional sublattice
$K_2$ of $\Gamma _i$ has minimum 8.

\begin{prop}
The vectors of norm $6$ in $\Gamma _i$ are of the form
$x+y$ with $x\in I_1 $ of norm $2$ and $y \in I_2$ of norm $4$.
\end{prop}

For all $i=1,\ldots ,9$
I computed representatives $v$ of
the orbits of $\Aut _{\Z[\alpha ]} (P_i) $
on the minimal vectors of $I_1$.
For each $v$ there is some $w\in I_2 $ such that
$v+w \in \Gamma _i$, moreover $w$ is unique modulo $K_2$.
So it remains to check that the minimum of the 48 dimensional lattice
$$ I(w):= \langle K_2 ,  w \rangle \leq I_2 $$
is $\geq 6$. This is done by enumerating all vectors of norm  4
in this lattice.

\begin{rem}
For $\Gamma _2 , \ldots , \Gamma _9$ the lattice
$I(w)$ contains vectors of norm $4$  for some $w$, summing
up to the number of vectors of norm $6$ in $\Gamma _i$ given
in Table 1.
Only for the lattice $\Gamma = \Gamma _1$
the representatives $w$ of all $15$ orbits of $\SL_2(25) $
provide lattices $I(w)$
of minimum $> 4$.
\end{rem}

\end{document}